\documentclass [12pt] {article}
\usepackage{amsfonts}
\usepackage{amssymb}

\newcommand{\qed} {\hspace {0.1in} \rule {1.5mm} {3.5mm}}

\newtheorem{lemma}{Lemma}[section]

\newtheorem{theorem}{Theorem}
\newtheorem{example}{Example}
\newtheorem{question}{Question}

\newtheorem{proposition}{Proposition}[section]
\newtheorem{definition}{Definition}[section]

\def\uv{\underline{v}}

\def\span{\mbox{span}_K\,}
\def\en{\mbox{End}\,}
\def\ke{\mbox{Ker}\,}
\def\ra{\mbox{Ran}\,}
\def\tr{\mbox{Tr}\,}
\def\coke{\mbox{Coker}\,}
\def\fred{\mbox{Fred}\,}

\def\dim{{\rm dim}}

\def\<{\langle}
\def\>{\rangle}

\def\noi{\vskip0.1in \noindent}
\def\proof{\smallskip\noindent{\it Proof.} }

\def\bZ{{\mathbb Z}}

\def\bN{{\mathbb N}}
\def\bC{{\mathbb C}}

\def\cA{\mbox{$\cal A$}}
\def\cB{\mbox{$\cal B$}}

\def\en{\mbox{End}}

\def\to{\rightarrow}

\begin{document}
\begin{center}
{.}
\vskip0.7in
{\Large {\bf The Ends of Algebras}}
\vskip0.3in
{\bf G\'abor Elek\footnote{The Alfred Renyi Mathematical Institute of
the Hungarian Academy of Sciences, P.O. Box 127, H-1364 Budapest, Hungary; Email:
elek@renyi.hu} and Aryeh Y. Samet-Vaillant\footnote{Department of Jewish Thought, 
The Hebrew
University of Jerusalem, Mount Scopus, 91905 Jerusalem, Israel;
Email: asamet@mscc.huji.ac.il}}
\vskip0.5in

{\large{\bf Abstract}} \end{center}\vskip0.1in
 We introduce the notion of ends for
algebras. The definition is analogous to the one in geometric group
 theory. We establish
some relations to growth conditions and cyclic cohomology.
\vskip 0.2in

\section{Introduction and preliminaries}
The goal of this paper is to introduce and study the notion of ends for
algebras. First let us recall the definition of ends of groups (see
e.g. Freudenthal (1945), Epstein (1962)).
Let $\Gamma$ be a group, then we call $P\subset \Gamma$ an almost invariant
set if for any $g\in \Gamma$, $Pg$ and $P$ differ only by finitely
many elements. That is $\{(Pg\cup P)\backslash P\}$ is
a finite set. Obviously, if $P$ is almost invariant, then $\Gamma\backslash P$
is almost invariant as well. If two almost invariant sets differ only by finitely many elements
then we call them equivalent, $P\simeq Q$. The equivalence classes of
almost invariant sets are the endsets.
We say that the group $\Gamma$ has $k$ ends if it can be partitioned into the
disjoint union of $k$ infinite almost invariant sets, but it can not be
partitioned into the disjoint union of $k+1$ infinite almost invariant sets.
If for any $k\geq 1$, the group has such partitions, then $\Gamma$ has
infinitely many ends. By the theorem of Hopf (1944), a group may have $0,1,2$ or
$\infty$ many ends. Finite groups have zero ends, the group $\bZ^2$ has
one end, the group $\bZ$ has two ends and the free group of two generators has infinitely
many ends. 

\noindent
Now let $K$ be a commutative field and $\cA$ be a unital $K$-algebra. We say that
a $K$-linear subspace $V\subset \cA$ is an {\bf almost invariant subspace} if for any $x\in\cA$,
$(Vx+V)/V$ is finite dimensional.
\begin{definition}
\label{d2}
The algebra $\cA$ has $k$ ends if $\cA=V_1\oplus V_2\oplus\dots\oplus V_k$,
where $\{V_i\}_{i=1}^k$ are infinite dimensional
almost invariant subspaces, but
$\cA$ can not be written as the direct sum of $k+1$ infinite dimensional
almost invariant subspaces. We say that $\cA$ has infinitely many ends, if
for any $k\geq 1$, $\cA$ has such direct sum decomposition. We apply the
following convention; finite dimensional algebras have zero ends.
\end{definition}
Let us see some examples.
\begin{example}
$\cA=K[x]$, the polynomial ring of one variable.
\end{example}
Let $S\subset \cA$ be an infinite dimensional almost invariant subspace.
Let $n>0$, such that $Sx+S\subset S+\span\{1,x,x^2,\dots,x^n\}$. Suppose
that $p_m\in S$ is a polynomial of degree $m>n$. Then $p_mx\in
 S+\span\{1,x,x^2,\dots,x^n\}$. That is $S$ contains a polynomial of degree $m+1$.
Inductively, $S$ contains a polynomial of degree $k$ for any $k\geq m$.
Hence $K[x]=S+\span\{1,x,x^2,\dots,x^{m-1}\}$. Therefore $K[x]$ has only one end.
\begin{example}
$\cA=K[x,x^{-1}]=K\bZ$.
\end{example}
If $S\subset \cA$ is an infinite dimensional almost invariant subspace,
then our previous example shows that $K[x]\cap S$ is either finite or
cofinite dimensional, similarly,  $K[x^{-1}]\cap S$ is either finite or
cofinite dimensional. Hence $K[x,x^{-1}]$ has two ends. 
\begin{example}
The $\ast$-algebra $\Sigma_n(K)=\oplus_{1\leq i \leq n} K[x_i]$
\end{example}
One can easily prove that if $\cA$ and $\cB$ are algebras, having respectively
$k$ and $l$
ends then $\cA\oplus\cB$ has $k+l$ ends.  Thus $\Sigma_n(K)$ has exactly $n$ ends.
\begin{example}
The algebra $\cA=K[x,y]/J$, where the ideal $J$ is generated by the monomials
of $x$ and $y$ containing at least $2$ $y$'s. 
\end{example}
Note that $S_k={\mbox span}\,\{x^k,x^ky, x^kyx, x^kyx^2,\dots\}$ is an almost
invariant subspace. Then
$$\cA=S_1\oplus S_2\oplus\dots S_n\oplus \{K\oplus S_{n+1}\oplus
S_{n+2}\dots\}$$
is a decomposition of $\cA$ into the direct sum of $n+1$ infinite dimensional
almost invariant subspaces. Hence $\cA$ has infinitely many ends.

\noindent

We shall define endclasses of algebras, similarly to the endsets of groups, and
identify the endclasses with the idempotents of a certain natural extension
of our algebras. Thus one-endedness of a certain algebra
is equivalent to the fact that this extension has no non-trivial idempotents.
Houghton (1972) proved that if a group is infinite then the number of ends
is exactly the dimension of the first cohomology group $H^1(\Gamma,\bC\Gamma)$
plus one. 
 Using a Fredholm-module construction 
we associate to the endclasses of $\cA$ elements of the first cyclic cohomology
group $HC^1(\cA)$. We also  prove that algebras of Gelfand-Kirillov
dimension $1$ have only finitely many ends.
Finally, we study the case of group algebras. 
\section{The Fredholm extension}
Let $K$ be a field and $K^\infty$ be the countable dimensional vector space
over $K$. Consider the endomorphism ring $\en_K(K^\infty)$. 
Let $I\subset \en_K(K^\infty)$ be the ideal of endomorphisms of finite rank.
If $T\in I$ then let $\widetilde{T}$ be the restriction of $T$ onto the
range of $T$. Then $\widetilde{T}$ is a finite dimensional linear
transformation hence $Trace({\widetilde{T}})$ is well-defined.
If $T'$ is the restriction of $T$ onto a finite dimensional
space containing the range of $T$, then of course $Trace(\widetilde{T})=
Trace (T')$. We shall denote $ Trace(\widetilde{T})$ by $\tr(T)$.
Obviously, $\tr: I\to K$ is a linear functional.
\begin{lemma}
\label{l4}
If $A\in\en_K(K^\infty)$ and $B\in I$, then $\tr(AB)=\tr(BA)$.
\end{lemma}
\proof Let $\{e_i\}_{i\geq 1}$ be a basis of $K^\infty$ such that
$e_1,e_2,\dots, e_n$ span the range of $B$ and the span of
$\{e_1,e_2,\dots,e_m\}$ contains the range of $AB$, where $m>n$.
Let $Be_i=\sum^\infty_{j=1} b_{ij} e_j$, $Ae_i=\sum^\infty_{j=1} a_{ij} e_j$.
Then $b_{ij}=0$ if $j>n$ and $a_{ij}=0$, if $j>m$ and $i\leq n$.
Hence $a_{ij} b_{ji}=0$, if either $i>m$ or $j>m$. Thus
$$
\tr(AB)=\sum^m_{i=1}\sum^\infty_{j=1} a_{ij}b_{ji}=
\sum^m_{i=1}\sum^m_{j=1} a_{ij}b_{ji}=
\sum^m_{i=1}\sum^\infty_{j=1} a_{ij}b_{ji}=\tr(BA)\quad\qed
$$
Now let us suppose that $\cA$ is a countable dimensional unital $K$-algebra.
For $x\in \cA$, denote by $l_x$ (resp. $r_x$) the left (resp. right)
multiplications by $x$.
Obviously, if $T\in\en_K(\cA)$ and $T$ commutes with $r_x$ for any $x\in\cA$,
then $T= l_{T(1)}$. We can identify the set $\{l_x\}_{x\in\cA}$ with $\cA$
itself and the set $\{r_x\}_{x\in\cA}$ with the opposite algebra $\cA^o$.
Let
$$
\widetilde{\cA}=\{T\in \en_K(\cA)\,\mid\, Tr_x-r_xT\in I,\,\mbox{if}\,
x\in\cA\}\, .
$$

Obviously $\cA\subset\widetilde{A}, I\subset\widetilde{A}$ and
$\widetilde{A}$ is a unital subalgebra of $\en_K(\cA)$.
\begin{definition}
The Fredholm extension of $\cA$ is defined as $\fred(\cA)=\widetilde{\cA}/I$.
\end{definition}
Note that
if for any $0\neq x\in\cA$, $l_x$ is of infinite dimensional range then
$\pi:\cA\to\fred(\cA)$ is injective. If $\cA$ is the polynomial ring $K[x]$,
then $\fred(\cA)$ is just the Toeplitz-algebra.

\begin{proposition}
\label{p6}
Suppose that $x$ is a Fredholm element,that is both $\ke (l_x)$ and
$\coke (l_x)$ is finite dimensional, then $[x]$ is invertible in
$\fred(\cA)$.
\end{proposition}
\proof
First of all note that if $W\subseteq\cA$ is a finite codimensional subspace
and $T\in\en_K(\cA)$, then there exists
a finite codimensional subspace $W_T\subset W$ such that $T(W_T)\subseteq W$.
Indeed, consider
the linear map $\pi_W\circ T: W\to \cA/W\,.$
Then $W_T=\ke(\pi_W\circ T)$ is clearly finite codimensional
and $T(W_T)\subseteq W$.
Now let $W^x\subset\cA$ be a subspace, $\ke (l_x)\oplus W^x=\cA$.
Then $ l_x\mid_{W^x}$ is injective. Consider
the finite codimensional subspace $l_x(W^x)$.
All elements in $l_x(W^x)$ can be uniquely written as
 $xw$,
where $w\in W^x$.
Let $\cA=l_x(W^x)\oplus B$, where $B$ is finite dimensional.
Define $T\in\en_K(\cA)$ to be $0$ on $B$ and $T(xw)=w$ on
$l_x(W^x)$.
Then $T\in\widetilde{\cA}$. Indeed, let $s\in\cA$ and let
$W^x_{r_s}\subset W^x$ be
a finite codimensional space as above, that is
 $r_s(W^x_{r_s})\subset W^x$. If
$w\in W^x_{r_s}$, then
$$Tr_s(xw)=T(xws)=ws\quad\mbox{and}\quad r_sT(xw)=ws\,.$$
Hence, $[T]$ is the inverse of $[x]$ in $\fred{\cA}$\qed

Let $\cA=V\oplus W$ be a decomposition of $\cA$ into the direct sum of
two almost invariant subspaces, where $V$ is infinite dimensional
. A decomposition
$\cA=V'\oplus W'$ is {\bf equivalent} to the decomposition above, if
$V$ and $V'$ respectively $W$ and $W'$ differ only by a finite dimensional
subspace. That is $(V+V')/V$, $(V+V')/V'$, $(W+W')/W$ and $(W+W')/W'$ are
all finite dimensional. We call the equivalence classes of such decompositions
{\bf the endclasses of the algebra $\cA$}. Note that $\cA=V\oplus W$ and
$\cA= W\oplus V$ are determining different endclasses. Of course, if $V$ is an
almost invariant subspace it might not have a complementary almost invariant
subspace at all. The following example shows that in some pathological
situations an almost invariant subspace can represent more than one endclasses.
\begin{example}
$\cA=\en_K(K^\infty)$.
\end{example}
Let $\{e_i\}_{i=1}^\infty$ be a basis of $K^\infty$.
Let $A=\span\{e_1, e_3, e_5, \dots\}, B=\span\{e_2, e_4, e_6, \dots\},
C=\span\{e_1+e_2, e_3+e_4, e_5+e_6,\dots\}$. For $X\subset K^\infty$,
let $$N_X=\{T\in \en_K{K^\infty}\,\mid\, T(v)=0\,\mbox{for any $v\in
X$}\}\,.$$
Then clearly, $N_A, N_B$ and $N_C$ are pairwise inequivalent almost
invariant subspaces, and $\en_K{K^\infty}=N_A\oplus N_B= N_A\oplus N_C=
N_B\oplus N_C\,.$

\begin{theorem}
The nonzero idempotents of $\fred(\cA)$ are in bijective correspondence
with the endclasses of $\cA$.
\end{theorem}
\proof
Let $\cA=V\oplus W$ be a decomposition of our algebra into the direct
sum of almost invariant subspaces. Let $\pi_V$ (resp. $\pi_W$) be the
projection
onto $V$ (resp. $W$). Then $\pi_V$ is an idempotent in $\cA$. Indeed
if $x\in \cA$, then by the observation in Proposition \ref{p6}, there exist
$R^V_x\subseteq V$, $R^W_x\subseteq W$ finite codimensional subspaces
such that $r_xR^V_x\subseteq V$ and $r_xR^W_x\subseteq W$.
Hence if $a\in R^V_x$, then $r_x\pi_V(a)=\pi_V(ax)=ax$,
$r_x\pi_W(a)=\pi_W(ax)=0$, thus $[\pi_V,r_x]\in I$.
If $A=V'\oplus W'$ an equivalent decomposition then clearly $[\pi_V]=
[\pi_{V'}]$ in $\fred(\cA)$.
Now, if $P\in \fred(\cA)$ is a non-zero idempotent, then let
$\widetilde{P}\in\widetilde{\cA}$ be a representative of $P$ and
$\widetilde{Q}\in\widetilde{\cA}$ be a representative of $1-P$.
It is easy to check that $\ra(\widetilde{P})$ and $\ra(\widetilde{Q})$
are almost invariant subspaces, such that $\ra(\widetilde{P})+
\ra(\widetilde{Q})$ is finite codimensional in $\cA$ and 
$\ra(\widetilde{P})\cap \ra(\widetilde{Q})$ is finite dimensional.
Hence there exist almost invariant subspaces $V,W$, such that
$\cA=V\oplus W$ and $V$ is equivalent to $\ra(\widetilde{P})$,
$W$ is equivalent to $\ra(\widetilde{P})$, where $V$ is infinite dimensional.
Then $[\pi_V]=P$, $[\pi_W]= 1-P$, thus the Theorem follows. \qed

\section{The noncommutative geometry of the ends}

In this section we apply the Fredholm module construction of Connes (1994). 
First
recall the notion of a graded differential algebra of length $1$.
Let $\cA$ be a unital $K$-algebra and $\Omega$ be an $\cA$-bimodule with
a $K$-linear map $\int:\Omega\to K$ and a $K$-homomorphism
$d:\cA\to \Omega$ satisfying the following conditions
\begin{itemize}
\item $d(ab)=a\cdot db+da\cdot b$, if $a,b\in\cA$.
\item $\int\, a\omega=\int\, \omega a$, if $a\in\cA$,\,$\omega\in \Omega$.
\item $\int\,da=0$, if $a\in\cA$.
\end{itemize}
Then $\tau:\cA\otimes_K\cA\to K$, $\tau(a,b)=\int a\cdot db$ defines a 
cyclic cohomology of $\cA$, $\tau\in HC^1(\cA)$. 
Indeed, one can easily check that 
$$\tau(a,b)=-\tau(b,a)$$
and
$$\tau(ab,c)-\tau(a,bc)+\tau(ca,b)=0\,.$$
Also, $\tau$ defines a character $\phi_\tau:GL_1(\cA)\to K$, by
 $\phi_\tau(a)=\tau(a^{-1},a)$.

\noindent
The classical example is of course the algebra of smooth complex valued 
functions on the unit circle. Then,
$$\tau(f,g)=\frac{1}{2\pi i}\int f(z) dg(z)\quad.$$
 $\phi_\tau(f)=\frac{1}{2\pi i}\int \frac{f'(z)}{f(z)} dz$
is just the winding number of $f$.

\noindent
Now let $\cA$ be a finitely generated $K$-algebra and let $\cA=V\oplus W$ be
a decomposition of $\cA$ into the direct sum of two infinite dimensional
almost invariant subspaces. Let $F$ be defined the following way,
$F(v)=v$, if $v\in V$, $F(w)=w$, if $w\in W$. That is $F\in\widetilde{\cA}$
is an endomorphism such that $F^2=1$.
Then $(\cA,I,F)$ is a Fredholm-module, where $I\in\en_K(\cA)$ is
the trace ideal. The coboundary operator $d:\cA\to I$ is given by
$dx=[F,r_x]$, $\int:I\to K$ is defined as $\int \omega=\tr(\omega)$.

\noindent
Clearly, $d(xy)=dx\cdot y+x\cdot dy$ and $\int a\omega=\int \omega a$.
One should only check that
$\int da=0$, if $a\in\cA$.
Pick a basis $\{e_i\}^\infty_{i=1}$ for $V$ and
a basis $\{f_i\}^\infty_{i=1}$ for $W$ such that the span
of $\{e_1,e_2,\dots,e_n,f_1,f_2,\dots, f_n\}$ contains the range
of $[F,r_a]$. Let $e_ia=v_i+w_i$, $f_ia=v'_i+w'_i$, where
$v_i,v'_i\in V$ and $w_i,w'_i\in W$. Then $[F,r_a]e_i=-2w_i$ and
$[F,r_a]f_i=2v_i$, hence $\tr [F,r_a]$ must be zero.
Therefore we associated a first cyclic cohomology class in $HC^1(\cA)$ to
the decomposition. Now if $\cA=V'\oplus W'$ is an equivalent decomposition
then $F-F'\in I$. Thus,

$$\int a (Fb-bF)-\int a(F'b-bF')=\int (a(F-F')b-ab(F-F'))=\int
(ab-ba)(F-F')\,.$$
This shows that $\tau-\tau'$ is a cyclic coboundary. Hence we proved the
following theorem.
\begin{theorem}
If $\{A,I,F\}$ is a triple as above associated to a decomposition $\cA=V\oplus
W$, then $p(a,b)=\tr (a[F,b])$ is a cyclic
cocycle. The cyclic cohomology element represented by $p$ in $HC^1(\cA)$
depends only on the endclass associated to the decomposition.
\end{theorem}

\noindent
Finally, let us show that at least in some cases one can obtain
non-trivial cyclic cohomologies.
Let $K[x,x^{-1}]=\cA=V\oplus W$, where $V=\span\{1,x,x^2,\dots\}$,
$W=\{x^{-1},x^{-2},\dots\}$.
Let us calculate $\phi_\tau(x)=\tr(x^{-1}[F,r_x])$.
Clearly $(Fr_x-r_xF) (x^k)=0$ if $k\neq -1$. On the other hand
$(Fr_x-r_xF) (x^{-1})=2$ thus $\tr(x^{-1}[F,r_x])=2.$ That is the 
associated cyclic cohomology is non-trivial.
\section{Algebras of linear growth}
In this section $\cA$ shall be a finitely generated $K$-algebra. Let $t_1,t_2,\dots t_n$
be a generator system of $\cA$ containing the unit. Let $S_1$ be the vector
space spanned by $\{t_1,t_2,\dots t_n\}$. Set $S_2=S_1S_1, S_3=S_2S_1$,
inductively $S_{n+1}=S_nS_1$. Then $S_1\subset S_2\subset\dots$ are finite
dimensional vector spaces and $\cup^\infty_{n=1}S_n=\cA$. Now we briefly
review some basic notions on the Gelfand-Kirillov dimension (see Krause and
Lenagan (2000)).
The growth function of $\cA$, relative to the generator system is the function
given by $f(n):=\dim (S_n)$, and the Gelfand-Kirillov dimension is
$$GK\dim(\cA)=\limsup\frac{\log(f(n))}{\log n}\,.$$
Note that the Gelfand-Kirillov dimension does not depend on the choice
of the generators. If $GK\dim(\cA)=1$, then by Bergman's theorem, $\cA$ in 
fact has linear growth, that is $f(n)\leq cn$, for some real number $c>0$.
Bergman also showed that the linear growth of an algebra
 is equivalent to the existence of a constant $K>0$, such that
$f(n+1)-f(n)<K$.

\noindent
By definition, if $V\subseteq A$ is an almost invariant subspace then there
exists $k_V>0$ such that $VS_1\subseteq V+S_{k_V}$.
\begin{lemma} 
\label{luj1}
Let $\cA=V\oplus W$ be a direct sum decomposition into almost invariant
subspaces, then for some $l>0$,
$S_n\subseteq (V\cap S_{n+l})\oplus (W\cap S_{n+l})$, holds for any $n>0$.
\end{lemma}
\proof
Let $k>0$ be an integer such that $VS_1\subseteq V+S_k$,
$WS_1\subseteq W+S_k$ and let $j>k$ be an integer such that 
$S_k\subseteq (V\cap S_j)\oplus (W\cap S_j)$. If $z\in S_k$,
$z=v+w, v\in V\cap S_j, w\in W\cap S_j$, then $zt_i=vt_i+wt_i$, where 
$vt_i\in (V\cap S_{j+1})+S_k$ and $t_iw\in (W\cap S_{j+1})+S_k$.
Hence $S_{k+1}\subseteq (V\cap S_{j+1})\oplus (V\cap S_{j+1})$. Inductively,
$S_{k+n}\subseteq (V\cap S_{j+n})\oplus (V\cap S_{j+n})$. \qed

\noindent
We call a subspace $V\subset \cA$ {\bf sparse} if for infinitely many integer
$m$, $V\cap S_m=V\cap S_{m+1}$.
\begin{lemma}
\label{luj2}
If $\cA$ is a finitely generated algebra and $V\subset\cA$ is a sparse almost
invariant subspace, then $V$ is finite dimensional. \end{lemma}
\proof
Let $\cA=V\oplus W$ as above , $m>j$ and $V\cap S_m=V\cap S_{m+1}$. Suppose
that $z\in S_k$, $z= v+ w$, $v\in V\cap S_m$, $w\in W\cap S_{m+1}$.
Then $zt_i\in (V\cap S_{m+1})\oplus (W\cap S_{m+2})=(V\cap S_m)\oplus
(W\cap S_{m+2}).$ Inductively, $S_n\subseteq (V\cap S_m)\oplus
(W\cap S_{m+n})$. Hence, $(V\cap S_m)\oplus W=\cA$, thus $V\subset
S_m$.\qed

\begin{theorem}
If $\cA$ is a finitely generated algebra of Gelfand-Kirillov dimension $1$, then
$\cA$ has only finitely many ends. Also if the growth function satisfies
$f(n)\leq Cn$, then $\cA$ has at most $C$ ends.
\end{theorem}
\proof
Suppose that for some $C>0$, $\dim_K S_n\leq C\,n$ and
 $\cA=V_1\oplus V_2\oplus \dots V_k$ and $\{V_i\}_{i=1}^k$ are infinite
dimensional almost invariant subspaces. Then none of the $V_i$'s are
sparse. Hence, for some $m>0$, $V_i\cap(S_{m+1}\backslash S_m)\neq \emptyset$,
for all $1\leq i \leq k$. Thus $\dim_K\,S_n\geq k(n-m)$, therefore
$k\leq C$. \qed

\noindent
Groups of linear growth must have two ends, however it follows from Example 3.
that the number of ends of algebras of linear growth may reach any positive
integer. Groups of non-exponential growth have only finitely many ends,
nevertheless the algebra of Example 4. in the Introduction, has quadratic growth
and infinitely many ends. 

\begin{question}
If $GK\dim(\cA)$ is finite and $\cA$ is finitely presented, is it possible
that $\cA$ has infinitely many ends ?
\end{question}

\section{The ends of group algebras}

Let $\Gamma$ be a group and $K$ is a commutative field.
If $V\subset\Gamma$ is an almost invariant subset then $KV\subset K\Gamma$
is an almost invariant subspace. Therefore if $\Gamma$ has at least $k$
ends, then $K\Gamma$ has at least $k$ ends as well. On the other hand let us
suppose that $\Gamma$ is a one-ended group, then $\Gamma\times C_2$ is
one-ended as well, where $C_2$ is the cyclic group of two elements. Then
$\frac{1+a} {2} $ is an idempotent in $\bC(\Gamma\times C_2)$, where $a$ is the
generator of $C_2$. In general, complex group algebras of infinite groups with torsion
always has more than one ends. Of course finite groups have zero ends and
its group algebras are finite dimensional, so they have zero ends as well.

\noindent
According to the Idempotent Conjecture, if $\Gamma$ is torsion-free, then
$\bC\Gamma$ has no non-trivial idempotent. The following question can be viewed
as a version of the Idempotent Conjecture.
\begin{question}
If  $\Gamma$ is torsion-free and has only one end, is it true that $\bC\Gamma$
has one end as well ?
\end{question}
Note that the answer is yes if and only if $\fred(\bC\Gamma)$ has no 
non-trivial idempotent. 
In the Abelian case we can answer the Question.

\begin{theorem}
If $\Gamma\simeq\bZ^n$, $n\geq 2$, then the group algebra $K\Gamma$
has only one end.
\end{theorem}
We introduce some notations. If $P\in K\bZ^n$ and $P=\sum_{\uv\in \bZ^n}k_{\uv}
\uv$, $k_{\uv}\in K$, where $[\uv]$ is the representative of the element
$\uv\in \bZ^n$ in $K\bZ^n$, thus $[\uv][\underline{w}]=[\underline{v}+
\underline{w}]$. Then let
$$T(P)=\{\uv\in \bZ^n\,\mid\,k_{\uv}\neq 0\}\quad.$$
If $\uv=(k_1,k_2,\dots,k_n)\in \bZ^n$, then $s(\uv)=\sup_{1\leq j \leq n} k_j$,
$i(\uv)=\inf_{1\leq j \leq n} k_j$. If $P\in K\bZ^n$, then
$s(P)=\sup_{\uv\in T(P)} s(\uv),$ \, $i(P)=\inf_{\uv\in T(P)} i(\uv).$
If $m\leq t$ are integers, then
$$A_{m,t}=\{(k_1,k_2,\dots,k_n)\in \bZ^n\,\mid\, m\leq k_j\leq t,\,
\mbox{if}\, 1\leq j \leq n\}\quad,\,A_m=\cup^\infty_{t=m} A_{m,t}\,.$$
The leading term of $P\in \bZ^n$, $l(P)$ is the greatest element of $T(P)$ in the
lexicographic ordering.
Furthermore, let $x_1=[(1,0,\dots,0)], x_2=[(0,1,0\dots,0)],\dots,
x_n=[(0,0,\dots,0,1)]\in K\bZ^n$.
\begin{lemma}
Let $V\subset K\bZ^n$ be an infinite dimensional almost invariant subspace.
Then if $d\in\bN$ is large enough, there exists $Q\in V$ such that
$l(Q)=(d,d,\dots,d)$.
\end{lemma}
\proof
By definition, if $V$ is an infinite dimensional almost invariant
subspace of $K\bZ^n$, then there exists a finite set
$\{\underline{b}_1, \underline{b}_2,\dots,\underline{b}_s\}$ such
that for any $R\in V, Rx_i= R'+L$, where
$R'\in V, L\in\span \{\underline{b}_1,
\underline{b}_2,\dots,\underline{b}_s\}$.
Let $i(b)=\inf_{1\leq j \leq s} i(\underline{b}_j),
s(b)=\sup_{1\leq j \leq s} s(\underline{b}_j).$ so the symmetric
difference of $T(R')$ and $T(R)+(0,0,\dots,1,\dots,0)$ is in the hypercube
$A_{i(b),s(b)}.$ Naively speaking we have control on the terms only
outside this hypercube.
Start with an element $P_1\in V$ with a term $\uv^1\notin A_{i(b),s(b)}\,.$
Then there exists $P_2\in V$ with a term $\uv^2\notin A_{i(b),s(b)}\,$ and
recursively we can obtain $P_l\in V$ with a term $\uv^l$, such that
$\uv^l>
\underline{b}_j$ in the lexicographic ordering, for any $1\leq j \leq s$.
Then repeating the process, we can obtain $Q\in V$, $l(Q)=(d,d,\dots,d)$,
$s(b)<d$.
\qed

\noindent
Now  fix an integer
$d>0$ such that for some $Q\in V$, $l(Q)=(d,d,\dots,d)$ and
$s(b)<d$. We also fix an integer
$m\leq 0$ such that $m<i(b)$ and
$m\leq\inf_{\uv\in T(Q)} i(\uv)$, for the same $Q\in V$.
Let $r=s(Q)-d$. Then $Qx_1=Q'+L'$, where $Q'\in V$ and $L'\in
\span \{\underline{b}_1,
\underline{b}_2,\dots,\underline{b}_s\}\,.$
 Obviously, $l(Q')=(d+1,d,\dots,d)$
and $s(Q')-(d+1)=r$. By induction, we can obtain the following lemma.
\begin{lemma}
For any $\uv\in\bZ^n$ such that $i(\uv)\geq d$, there exists $Q\in V$
such that $l(Q)=\uv$, $T(Q)\subset A_m$ and $s(Q)-s(\uv)=r$. 
\end{lemma}
\begin{lemma}
$$\lim_{t\to\infty}
\frac{\dim_K(V\cap \span(A_{m,t}))}{|A_{m,t}|}=1$$
\end{lemma}
\proof
First note that if $P_1,P_2,\dots,P_l\subset K\bZ^n$ have different
leading terms then they form an independent system over $K$. By the previous
lemma, if $t$ is large enough then for any $\uv$, $d\leq i(\uv)\leq t-r$,
there exists $Q\in V$ such that $l(Q)=\uv$ and $T(Q)\subset A_{m,t}$.
Hence,
$$\frac{\dim_K(V\cap \span(A_{m,t}))}{|A_{m,t}|}\geq
\frac{|A_d\cap A_{m,t-r}|}{|A_{m,t}|}\,.$$
However, $\lim_{t\to\infty} \frac{|A_d\cap A_{m,t-r}|}{|A_{m,t}|}=1$\qed

\noindent
Now the Theorem follows. Indeed, by the previous lemma, any two infinite
dimensional almost invariant subspaces have non-zero intersections. \qed

\section{Some questions on domains}

Among the intensively studied ``geometric'' properties of groups, one
has the growth, amenability and ends (see Gromov (1993)).
As mentioned in the course of this paper, these concepts are
closely related. The three concepts have been introduced in the
context of algebras (See Elek (2003) and
Samet (2000) for discussions on the amenability of algebras).
Groups of polynomial growth has
integer degree of growth. Algebras of polynomial growth may have 
any real number greater or equal two as Gelfand-Kirillov dimension
(Krause-Lenagan (2000)).
 It is an old problem, whether all
 domains of polynomial growth have integer Gelfand-Kirillov
dimension.

\noindent
By Hopf Theorem, a group may have zero, one, two or infinitely many ends.
All the domains we can check has zero, one, two or
infinitely many ends. So we can ask a similar question.
\begin{question}
Is it true, that a domain may have only zero, one, two or infinitely many
\\ ends ?
\end{question}

\noindent
Finally, it is well known, that amenable groups have finitely many
ends. The algebra in Example 3. has quadratic growth and so it is amenable.
On the other hand, it has infinitely many ends. Again, we can ask about
the similarity of domains to groups.

\begin{question}
Is it true, that any amenable domain has only finitely many ends ?
\end{question}

\vskip0.2in
\noindent
{\large{\bf References}}
\noi
Connes, A. (1994). {\it Noncommutative Geometry.}
Academic Press Inc: San Diego.
\noi
Elek, G. (2003). {\it The amenability of affine algebras.}
J. Algebra 264:469-478.
\noi
Epstein, D. (1962). {\it Ends}.
In: {\it Topology of 3-manifolds and related results},
Ford M. Ed., Prentice Hall, pp. 110-117.
\noi
Freudenthal, H. (1945). {\it \"Uber die Enden diskreter R\"aume und Gruppen}
Comment. Math. Helv. 17:1-38.
\noi
Gromov, M. (1993). {\it Asymptotic invariants of infinite groups.}
in: {\it Geometric Group Theory}, Graham A. Niblo and Martin A. Roller,
Eds., Vol. 2, London Mathematical Society Lecture Note Series 182,
Cambridge Univ. Press, Cambridge, UK, pp. 1-295.
\noi
Hopf, H. (1944). {\it Enden offener R\"aume und unendliche diskontinuierliche
Gruppen.} Comment. Math. Helv. 16:81-100.
\noi
Houghton, C. (1972). {\it Ends of groups and the associated first cohomology
group.} J. London Math. Soc. 6:81-92.
\noi
Krause, G., Lenagan, T.(2000) {\it Growth of algebras and the Gelfand-Kirillov
dimension} (Revised edition), Graduate Studies in Mathematics, 22,
American Mathematical Society: Providence, RI.
\noi
Samet-Vaillant, A. Y. (2000) {\it $C^*$-algebras, Gelfand-Kirillov dimension
and Folner sets.} J. Funct. Anal. 171:346-365.

\end{document}